\documentclass[11pt]{article}
\usepackage{amsmath,epsfig, amsthm, amssymb, enumerate,algorithmic,algorithm}
\usepackage{fullpage}
\usepackage{graphicx}
\usepackage{subfigure}
\usepackage{amsmath, amssymb, epsfig}
\usepackage{psfrag}
\usepackage{bm}
\usepackage{mathrsfs}
\usepackage{color}

\DeclareMathOperator*{\argmax}{arg max}

\newcommand{\new}[1]{{ #1}}

\newcommand{\R}{\mathbb{R}}

\renewcommand{\P}{\mathbb{P}}

\newcommand{\lp}[1]{\langle}
\newcommand{\rp}[1]{\rangle}

\newcommand{\istar}{i^*}
\newcommand{\xtilde}{\tilde{x}}
\newcommand{\atilde}{\tilde{a}}
\newcommand{\Vtilde}{\tilde{V}}
\newcommand{\ctilde}{\tilde{c}}
\newcommand{\Sc}{S^{c}}
\newcommand{\Vts}{\Vtilde_{S}}
\newcommand{\Vtsc}{\Vtilde_{S^c}}

\newcommand{\vs}{v_s}

\newcommand{\vtilde}{\tilde{v}}
\newcommand{\xsharp}{x^\#}
\newcommand{\zsharp}{z^\#}
\newcommand{\xtildesharp}{\xtilde^\#}

\newtheorem{theorem}{Theorem}[]
\newtheorem{lemma}[theorem]{Lemma}

\newcommand{\setdoublespace}{} 

\begin{document}
\title{Scaling Law for Recovering the Sparsest Element in a Subspace}

\author{Laurent Demanet and Paul Hand \\ $\;$ \\ Massachusetts Institute of Technology, Department of Mathematics, \\ 77 Massachusetts Avenue, Cambridge, MA 02139}
\date{October 2013, Revised May 2014}		

\maketitle

\begin{abstract}
{
\setdoublespace
We address the problem of recovering a sparse $n$-vector within a given subspace.  \new{This problem is a subtask of some approaches to dictionary learning and sparse principal component analysis.  Hence, if we can prove scaling laws for recovery of sparse vectors, it will be easier to derive and prove recovery results in these applications.  In this paper, we present a scaling law for recovering the sparse vector from a  subspace that is spanned by the sparse vector and $k$ random vectors.  }
We prove that the sparse vector will be the output to one of $n$ linear programs with high probability if its support size $s$ satisfies $s \lesssim n/\sqrt{k \log n}$. The scaling law still holds when the desired vector is approximately sparse. 
To get a single estimate for the sparse vector from the $n$ linear programs, we must select which output is the sparsest.  This selection process can be based on any proxy for sparsity, and the specific proxy has the potential to improve or worsen the scaling law.  If sparsity is interpreted in an $\ell_1/\ell_\infty$ sense, then the scaling law can not be better than $s \lesssim n/\sqrt{k}$.  Computer simulations show that selecting the sparsest output in the $\ell_1/\ell_2$ or thresholded-$\ell_0$ senses can lead to a larger parameter range for successful recovery than that given by the $\ell_1/\ell_\infty$ sense.  
}
{sparsity, linear programming, signal recovery, sparse principal component analysis, dictionary learning.}
\end{abstract}

\section{Introduction} 
\new{
We consider the task of finding the sparsest nonzero element in a subspace.  That is, given a basis for the subspace $W \subset \R^n$, we seek the minimizer of the problem
\begin{align}
\min \| z\|_0 \text{ s.t. } z \in W, z \neq 0, \label{l0-min}
\end{align}
where $\|z\|_0$ is the number of nonzeros in the coefficients of $z$.  
Because problem \eqref{l0-min} is NP-hard over an arbitrary subspace \cite{CP1986}, we study a convex relaxation over particular subspace models for which we can prove recovery results.

This problem is a special case of the more general task of finding many independent sparse vectors in a subspace \cite{GN2010}.   This general task has many applications, of which some will require a full basis of sparse vectors and others will require only a few sparse vectors.
}

\subsection{Applications}

\new{
Finding sparse vectors in a subspace is an important subtask in several sparse structure recovery problems.  Two examples are dictionary learning and sparse principal component analysis (PCA).
\begin{itemize}
\item Dictionary Learning with a square, sparsely used dictionary \cite{SWW2012}.  Consider an $m \times n$ matrix $Y$, where $n > m$ and each column of $Y$ can be written as a sparse linear combination of $m$ unknown and independent dictionary elements.  That is to say, $Y = A X$, where $A$ is invertible, the columns of $A$ are the dictionary elements, and the rows of $X$ are sparse.   The goal is to find $A$ and $X$ from only the knowledge of $Y$.  In order to find this decomposition, we note that the row span of $Y$ is the same as that of  $X$.  Hence, a process that finds a sparse basis of a subspace could be used to find $X$ from the row space of $Y$.  The dictionary $A$ can then be directly computed.  Such a process based on linear programming was used in \cite{SWW2012} to study this dictionary learning problem.  
\item Sparse PCA in the infinite data limit.  Consider the following noisy data:
\begin{align}
x_i \sim N(0, \beta P_W + I), \text{ for } i = 1 ... m, \label{sparse-pca-data}
\end{align}
where $\beta$ is a signal-to-noise ratio and $P_W$ is the projection matrix onto a subspace $W$.
The task of PCA is to identify the directions along which there is the most variability in the data.  For \eqref{sparse-pca-data}, these directions are the subspace $W$, which we interpret as the signal.  In applications involving gene expression data, these directions may consist of the combinations of genes that distinguish different types of cancer cells \cite{LdA2010}, for example.  In some applications, it is desirable that the discovered signal involve few variables \cite{ZAEG2010}.  Hence, sparse signal vectors are of particular interest.  For example, this sparsity allows easier interpretation of the relevant genes for the cancer discrimination task.  It is thus useful to find the sparse signal directions within $W$ based on the measurements $x_i$, which is a sparse PCA problem \cite{ZHT2006,DGJL2007, LT2013}.  This problem is most interesting for a small number of samples, but it is nontrivial even in the infinite data limit.  In this limit, $P_W$ is measured exactly, and sparse PCA reduces to finding the sparsest element in the subspace $W$.
\end{itemize}

To prove rigorous recovery results for problems like these, it is important to have scaling laws for when the sparsest element subtask succeeds.

}

\subsection{Recovery by $n$ linear programs}

\new {

The difficulty of finding sparse elements in a subspace stems from the nonconvexity of both the $\ell_0$ objective and the nonzero constraint in \eqref{l0-min}.  As is standard in sparse recovery,  we consider the  $\ell_1$ relaxation of the $\ell_0$ sparsity objective.  The nonzero constraint must now be replaced with some proper normalization in order to prevent outcomes arbitrarily close to zero.  

A natural normalization is to search for an element with unit Euclidean length.  Unfortunately, the resulting problem is nonconvex, and methods for solving it directly get stuck in local minima \cite{ZP2001}.  This approach may be convexified; for more details, see Section \ref{sec:discussion} of the present paper.

A convex normalization is to search for an element where  a particular coefficient is set to unity.  Such an approach should work best if we could normalize the largest coefficient to be $1$.  Unfortunately, we do not know which component corresponds to the largest coefficient of the sparsest element.  Hence, we attempt to recover the sparse element by separately normalizing each of the $n$ components.  That is, we attempt to find the sparsest element in $W$ by solving the collection of $n$ linear programs
\begin{align}
\min \|z\|_1 \text{ such that } z \in W, z(i) = 1, \label{optimization-problem-z-in-v}
\end{align} 
where $1 \leq i \leq n$. This is the approach introduced by Spielman et al. in \cite{SWW2012}.
}
In order to get a single estimate from the optimizers of these $n$ linear programs, we need a selector that returns the output that is the `sparsest.'  Here, the sparsity of a vector $z$ may be interpreted in one many precise senses, \new{such as the strict-$\ell_0$ sense, $\|z\|_0$; a thresholded-$\ell_0$ sense, $\#\{i \mid |z(i)| \geq \epsilon \}$; the $\ell_1/\ell_\infty$ sense,  $s = \|z\|_1/\|z\|_\infty$; or the $\ell_1/\ell_2$ sense, $ \|z\|_1/\|z\|_2$.  In principle, this interpretation of sparsity may alter the scaling law of the overall recovery process.  }

\subsection{Subspace Models}
\new{
Because finding the sparsest element in an arbitrary space is NP-hard, we restrict our attention to certain subspace models.  Two such models are the Bernoulli-Gaussian model and a planted-random model.
\begin{itemize}
\item Bernoulli-Gaussian model \cite{SWW2012}.  Consider a collection of  vectors where each entry is nonzero with a fixed probability.  Suppose the values of these nonzero coefficients are given by independent normal random variables.  The Bernoulli-Gaussian subspace is given by the span of these vectors.  If the subspace dimension and the expected sparsity of each vector is small enough, then the sparsest basis provides exactly these vectors, up to scaling.  Spielman et al. \cite{SWW2012} provide a scaling law for successful recovery with the linear programs \eqref{optimization-problem-z-in-v} under this subspace model.
\item Planted-random model.  Consider a sparse vector and a collection of independent Gaussian random vectors.  The planted-random subspace is given by the span of all these vectors.  If the sparse vector is sparse enough, then it is the unique sparsest element with probability 1, up to scaling.  To the best of our knowledge, this model has not been studied in the literature.

\end{itemize}

In the present paper we will provide scaling laws for successful recovery under the planted-random subspace model.  Concisely, the problem is as follows:
\begin{align}
\text{Let}&:  v \in \R^n, \vtilde_{j} \sim N(0, I_n) \text{ for } 1 \leq j \leq n \notag \\
\text{Given}&: \{w_1, \ldots, w_{k+1}\} \text{ a basis for } W = \text{span} \{ v, \vtilde_{1}, \ldots, \vtilde_{k} \} \label{planted-random-problem}\\ 
\text{Find}&: v \notag
\end{align}
The vector $v$ is the unique sparsest element in the space with probability 1, up to scaling, provided that $\|v\|_0 < n -k$.  
}
\subsection{Exact Recovery by a Single Program} \label{sec:exact-recovery}

Solving \eqref{planted-random-problem}  by the $n$ linear programs \eqref{optimization-problem-z-in-v} requires that at least one of those programs recovers $v$.  Roughly speaking, the program of form \eqref{optimization-problem-z-in-v}  that is most likely to succeed corresponds to an $\istar \in \argmax_i |v(i)|$.  This is because $v/v(i)$ is feasible for \eqref{optimization-problem-z-in-v} and has the smallest $\ell_1$ norm when $i =\istar$.  
The following theorem provides a scaling law under which $v$ is exactly recovered by this instance of  \eqref{optimization-problem-z-in-v}.

\begin{theorem} \label{thm-recovery-exact}
Let $k \leq n/32.$  There exists a universal constant $c$ such that for sufficiently large n, 
\begin{align}
\|v \|_0 \leq c \frac{n / \sqrt{\log n}}{\sqrt{k}} \Rightarrow \frac{ v}{v(\istar)} \text{ is the unique solution to \eqref{optimization-problem-z-in-v} for } i = \istar, \label{thm-exact-statement}
\end{align}
with probability at least $1 - \new{\tilde{\gamma}_1 e^{-\tilde{\gamma}_2 n /2}} - k e^{-\lfloor c \sqrt{n/\log n}\rfloor} - \frac{k}{n^2}$.  \new{Here, $\tilde{\gamma}_1$ and $\tilde{\gamma}_2$ are universal constants.  }


\end{theorem}


From the scaling law, we observe the following scaling limits on the permissible sparsity in terms of the dimensionality of the search space: \begin{alignat}{3}
k \text{ on the order of } &1 &&\Longrightarrow \|v\|_0 \lesssim n / \sqrt{\log n} \label{scaling-k-first}\\
k \text{ on the order of } &n &&\Longrightarrow \|v\|_0  \lesssim \sqrt{n} / \sqrt{\log n} \label{scaling-k-last}
\end{alignat}
That is, a search space of constant size permits the discovery of a vector whose support size is almost a constant fraction of $n$.  Similarly, a search space of fixed and sufficiently small fraction of the ambient dimension allows recovery of a vector whose support size is almost on the the order of the square root of that dimension.  

\new{
Roughly speaking, sparse recovery succeeds because the planted vector is smaller in $\ell_1$ than any item in the random part of the subspace.  The full minimization problem then reduces to computing the minimal value of \eqref{optimization-problem-z-in-v} for a random subspace with no planted vector.  In more precise terms, recall that $|v(\istar)| = \| v \|_\infty$ and $v/v(\istar)$ is feasible for \eqref{optimization-problem-z-in-v} with $i = \istar$.  A necessary condition for successful recovery is that the sparse vector is smaller in $\ell_1$ than the minimum value attainable by the span of the random vectors:
\begin{align}
\frac{\|v\|_1}{\|v\|_\infty} \leq \min \|z\|_1 \text{ such that } z \in \text{span} \{\vtilde_{1}, \ldots, \vtilde_{k}\}, z(\istar) = 1,   \label{rough-condition-recovery}
\end{align}
We will show in Section \ref{sec-thm-gist} that the right hand side of \eqref{rough-condition-recovery} scales like  $n/\sqrt{k}$ when $k$ is at most some constant fraction of $n$.  As $\|v\|_1/\|v\|_\infty \leq \|v\|_0$ for all $v$, and the equality is attained for some $v$, we conclude that high probability recovery of arbitrary $v$ is possible only if $\|v\|_0 \lesssim n/ \sqrt{k}$.

Because of this necessary condition \eqref{rough-condition-recovery}, the scaling law between $n, k,$ and $\|v\|_0$ in Theorem \ref{thm-recovery-exact} can not be improved, except for the logarithmic factor.  The scaling could conceivably be improved in the related context where we only seek some $i$ for which $v/v(i)$ is the unique solution to \eqref{optimization-problem-z-in-v}.  

}

\subsection{Stable Recovery by a Single Program}

The linear programs \eqref{optimization-problem-z-in-v} can also recover an approximately sparse $v$. That is, if $v$ is sufficiently close to a vector that is sparse enough in an $\ell_0$ sense, we expect to recover something close to $v$ by solving \eqref{optimization-problem-z-in-v} with $i = \istar \in \argmax_i |v(i)|$.  Let $\vs$ be the best $s$-sparse approximation of $v$.  The following theorem provides a scaling law under which $v$ is approximately the output to the $i = \istar$ instance of the program \eqref{optimization-problem-z-in-v}.  


\begin{theorem} \label{thm-recovery-noisy}
Let $k \leq n/32.$  There exists universal constants $c,C$ such that for sufficiently large n,  for $s = \lfloor c \frac{n / \sqrt{\log n}}{\sqrt{k}}\rfloor$, and for $i=\istar$, any minimizer $\zsharp$ of \eqref{optimization-problem-z-in-v} satisfies
\begin{align}
\left \|\zsharp - \frac{v}{v(\istar)} \right \|_2 \leq C \frac{\sqrt{k \log n}}{\sqrt{n}} \frac{\|v - \vs \|_1}{\|v \|_\infty} \label{final-error-estimate}
\end{align}
with probability at least $1  - \new{ \tilde{\gamma}_1 e^{-\tilde{\gamma}_2 n/2} } - k e^{-\lfloor c \sqrt{n / \log n}\rfloor} - k/n^2$. 
\end{theorem}

The dependence on $k$ and $n$ in \eqref{final-error-estimate} is favorable provided that $k \lesssim n / \log n$.  In the case that $k \sim n$, the error bound has a mildly unfavorable constant, growing like $\sqrt{ \log n}$.   The \new{$\sqrt{k/n}$} behavior of the error constant plays the roll of the $1/\sqrt{s}$ term that arises in the noisy compressed sensing problem \cite{CRT2005}.  The estimate \eqref{final-error-estimate} is slightly worse, as $\sqrt{k/n} \sim  k^{1/4}/\sqrt{s}$, ignoring logarithmic factors.  We believe that the $k$ and $n$ dependence of the error bound could be improved.  

\subsection{\new{Selecting the Sparsest Output}} \label{sec:sparsest-output}

Successful recovery of $v$ by solving the $n$ linear programs \eqref{optimization-problem-z-in-v} requires both that $v$ is the output to $\eqref{optimization-problem-z-in-v}$ for some $i$, and that $v$ is selected as the `sparsest' output among all $n$ linear programs.  We now comment on effects of selecting the sparsest output under several different interpretations of sparsity: $\ell_1/\ell_\infty$, $\ell_1/\ell_2$, and thresholded-$\ell_0$.  Computer simulations in Section \ref{sec:simulations} show that the precise sense in which sparsity is interpreted can substantially affect recovery performance, especially when the planted-random subspace has small dimension.  The worst empirical performance is exhibited by $\ell_1/\ell_\infty$, and the best is by thresholded-$\ell_0$, though this method has the drawback of requiring a threshold parameter.

\subsection{Discussion} \label{sec:discussion}
\new {
We now compare several approaches for solving the sparsest element problem \eqref{planted-random-problem}.   The main difference among these approaches is the way they relax sparsity into something tractable for optimization problems.  We then compare the sparsest element problem to compressed sensing.  
}

\new{
As the $\ell_1/\ell_\infty$ ratio is a proxy for sparsity,} it is natural to try to find the sparsest nonzero element in the subspace $W$ by solving

\begin{align}
\min \frac{\| z\|_1}{\|z\|_\infty} \text{ such that } z \in W, z \neq 0.  \label{l1linf-prob}
\end{align}
\new{
While this problem is not convex, it is the same as solving the $n$ programs \eqref{optimization-problem-z-in-v} and selecting the smallest output in an $\ell_1$ or an $\ell_1/\ell_\infty$ sense.  Geometrically, this approach corresponds to replacing an $\ell_\infty$ constraint ball by separate hyperplanes along each of its faces.
}

We note that the scaling for successful recovery with \eqref{l1linf-prob} can be no better than than $s \lesssim n/\sqrt{k}$.  To see this, observe that the right hand side of \eqref{rough-condition-recovery} provides an upper bound on the minimal value of $\|z\|_1 / \|z\|_\infty$ for $z \in \text{span} \{\vtilde_{1}, \ldots, \vtilde_{k}\}.$   Hence, outside the scaling $s \lesssim n/\sqrt{k}$, there would be random vectors that are sparser than $v$ in the $\ell_1/\ell_\infty$ sense.
 We do not recommend simply solving  \eqref{l1linf-prob} as an approach for finding the sparsest vector in a subspace because computer simulations reveal that recovery can be improved merely by changing the sparsity selector to $\ell_1/\ell_2$ or thresholded-$\ell_0$.

\new{
As the $\ell_1/\ell_2$ norm ratio is also a proxy for sparsity, it is also natural to try to solve the sparsest element problem by optimizing
\begin{align}
\min \frac{\| z\|_1}{\|z\|_2} \text{ such that } z \in W, z \neq 0.  \label{l1l2-prob}
\end{align}
This problem is also nonconvex, though it can be convexified by a lifting procedure similar to \cite{CESV2011}.  After such a procedure, it becomes an $n \times n$ semidefinite program.  While the procedure squares the dimensionality, it may give rise to provable recovery guarantees under a scaling law.    We are unaware of any such results in the literature.  In a sense, the present paper can be viewed as a simplification of this $n\times n$ matrix recovery problem into $n$ linear programs on vectors, provided that we are willing to change the precise proxy for sparsity that we are optimizing.  

A more elementary approach to solving the sparsest element problem \eqref{planted-random-problem} is to simply choose the largest diagonal elements of the projector matrix onto $W$.  If the sparse element $v$ is sparse enough, then this process will recover its support with high probability.  While this method may have a similar scaling as that in Theorem \ref{thm-recovery-exact}, we anticipate that it is less robust because it capitalizes directly on specific properties of the distribution of planted-random subspaces.  As an analogy, consider the problem of identifying a planted clique in a random graph.  A thresholding based algorithm \cite{K1995} performs at the best known tractable scaling, up to logarithmic factors; however, its performance is much worse with semirandom models \cite{FK2000}, where an adversary is able to modify the planted random graph within some constraints.  For the sparsest element problem, we leave a detailed study of the thresholding approach to future research.  
}

\new{
We now compare the sparsest element problem to compressed sensing.  The linear programs \eqref{optimization-problem-z-in-v} could be viewed as separate compressed sensing problems.  From this perspective, the sensing matrix of each would have many rows that are orthogonal to the planted sparse vector, hence complicating an analysis based on restricted isometries.  Because of this difficulty, we do not use a compressed sensing perspective for the proofs of the recovery theorems.  
}

\new{The sparsest element problem results in different qualitative scalings than those of compressed sensing.  For example, if sparsity is minimized in the $\ell_1/\ell_\infty$ sense, the scaling law can not be better than $s \lesssim n/\sqrt{k}$.  This scaling depends strongly on $k$.  From a compressed sensing perspective, one might guess that the recoverable sparsity is on the order of the number of measurements.  Because the problems \eqref{optimization-problem-z-in-v} can be written with $n-k+1$ equality constraints, this guess would provide the scaling $s \lesssim n$, which is impossible when minimizing sparsity in an $\ell_1/\ell_\infty$ sense.
 }
 


\subsection{\new{Notation}}

\new{
Let $W$ be the planted-random subspace spanned by $v$ and $\vtilde_1, \ldots, \vtilde_{k}$, as given in \eqref{planted-random-problem}. Let the matrix $V = [v, \vtilde_{1}, \vtilde_{2}, \ldots, \vtilde_{k}]$ have these vectors as columns.  Let $\Vtilde = [ \vtilde_{1}, \vtilde_{2}, \ldots, \vtilde_{k}]$, which gives $V = [v, \Vtilde]$.  
Let $v(i)$ be the $i$th component of the vector $v$.  
Let $V_{\istar,:}$ be the $\istar$-th row of $V$.  Write $V_{\istar,:} = [1, \atilde^t]$, where $\atilde \in \R^k$.   For a set $S$, write $\Vts$ and $\Vtsc$ as the restrictions of $\Vtilde$ to the rows given by $S$ and $\Sc$, respectively.  Write $x \new{\asymp} y$ when there exists positive $c$ and $C$, which are independent of $n$ and $k$,  such that $cy \leq x \leq Cy$. 
}

\section{Proofs} \label{sec:this}
To prove the theorems, we note that \eqref{optimization-problem-z-in-v}, \eqref{thm-exact-statement}, \eqref{final-error-estimate}, and the value of $\istar$ are all invariant to any rescaling of $v$.  Without loss of generality, it suffices to take $\|v\|_\infty = 1$ and $v(\istar) = 1$.
Our aim is to prove that $v$ is or is near the solution to \eqref{optimization-problem-z-in-v} when $i = \istar$.  We begin by noting that $W = \text{range}(V)$.  Hence, changing variables by $z = Vx$, \eqref{optimization-problem-z-in-v} is equivalent to 
\begin{align}
\min \|Vx\|_1 \text{ such that } Vx (\istar)= 1 \label{optimization-problem-vx}
\end{align}
when $i = \istar$.  Note that $Vx(\istar)$ refers to the $\istar$ component of Vx.  We will show that $x=e_1$ is the solution to \eqref{optimization-problem-vx} in the exact case and is near the solution in the noisy case.  Write $x = [x(1), \xtilde]$ in order to separately study the behavior of $x$ on and away from the first coefficient.  Our overall proof approach is to show that if $n$ is larger than the given scaling, a nonzero $\xtilde$ gives rise to a large contribution to the $\ell_1$ norm of $Vx$ from coefficients off the support of $v$.

\subsection{Scaling Without a Sparse Vector} \label{sec-thm-gist}

In this section, we derive the scaling law for the minimal $\ell_1$ norm attainable in a random subspace when one of the components is set to unity.  This law is the key part of the justification that the scaling in the theorems can not be improved except for the logarithmic factor.  As per Section \ref{sec:discussion}, it also provides the proof of the best possible scaling when minimizing sparsity in an $\ell_1/\ell_\infty$ sense.  We also present the derivation for pedagogical purposes, as it contains the key ideas and probabilistic tools we will use when proving the theorems.  The rest of this section will prove the following lemma.  

\begin{lemma}  Let $\Vtilde$ be an $n \times k$ matrix with i.i.d. ${N}(0,1)$ entries, where $k \leq n / 16$.  With high probability, 
\begin{align}
\frac{n}{\sqrt{k}} \new{\asymp} \min \| \Vtilde \xtilde \|_1 \text{ such that } \Vtilde \xtilde (\istar) = 1.  \label{scaling-law-component-normalization}
\end{align}
The failure probability is exponentially small in $n$ and $k$.  
\end{lemma}



\begin{proof}
Because $\text{range}(\Vtilde)$ is a $k$-dimensional random subspace, we can appeal to the uniform equivalence of the $\ell_1$ and $\ell_2$ norms, as given by the following lemma.

\begin{lemma} \label{lemma-l1-l2-random-subspace} Fix $\eta < 1$.  For every $y$ in a randomly chosen (with respect to the natural Grassmannian measure)  $\eta n$-dimensional subspace of $\R^n$, $$c_\eta \sqrt{n} \|y\|_2 \leq \|y\|_1 \leq \sqrt{n} \|y\|_2$$ with probability $1 - \gamma_1 e^{-\gamma_2 n}$ for universal constants $\gamma_1, \gamma_2$. 
\end{lemma}
\noindent This result is well known \cite{FLV1977, K1977}.  Related results with different types of random subspaces can be found at \cite{LS2007, AM2006, L2008,GLR2010}.  Thus, with high probability,
\begin{align}
\|\Vtilde \xtilde\|_1 \new{\asymp} \sqrt{n} \|\Vtilde \xtilde \|_2 \text{  for all }\xtilde. \label{l1-l2-equiv}
\end{align}

We now appeal to nonasymptotic estimates of the singular values of $\Vtilde$.   
Corollary 5.35 in \cite{V2010} gives that for a matrix $A \in \R^{n \times k}$ with $k\leq n/16$ and i.i.d. $N(0,1)$ entries,
\begin{align}
\P\Bigl(\frac{\sqrt{n}}{2}  \leq \sigma_\text{min}(A) \leq \sigma_\text{max}(A) \leq \frac{3 \sqrt{n} }{2}\Bigr) &\geq 1 -  2 e^{-n/32}. \label{sigma-min-random-matrix}
\end{align}
Thus, with high probability,
\begin{align}
\|\Vtilde \xtilde\|_2 \new{\asymp} \sqrt{n} \|\xtilde\|_2. \label{sing-vals-control}
\end{align}
Combining \eqref{l1-l2-equiv} and \eqref{sing-vals-control}, we get $\|\Vtilde \xtilde\|_1 \new{\asymp} n \| \xtilde\|_2$ with high probability.  Hence, the minimum values of the following two programs are within fixed constant multiples of each other:   
\begin{align}
\min \| \Vtilde \xtilde \|_1 \text{ such that } \Vtilde_{\istar, :} \xtilde = 1 \quad &\new{\asymp} \quad
\min n \| \xtilde \|_2 \text{ such that } \Vtilde_{\istar, :} \xtilde = 1. \label{two-progs}
\end{align}
By the Cauchy-Schwarz inequality and concentration estimates of the length of a Gaussian vector,  any feasible point in the programs \eqref{two-progs} satisfies
\begin{align}
\|\xtilde\|_2 \geq \frac{1}{\|\Vtilde_{\istar,:}\|_2} \new{\asymp} \frac{1}{\sqrt{k}}, \label{cauchy-schwarz-length-scaling}
\end{align}
with failure probability that decays exponentially in $k$.  Considering the $\xtilde$ for which the inequality in \eqref{cauchy-schwarz-length-scaling} is achieved, we get that the minimal value to the right hand program in \eqref{two-progs} scales like $n/\sqrt{k}$, proving the lemma.
\end{proof}

\subsection{Proof of Theorem \ref{thm-recovery-exact}}

The proof of Theorem \ref{thm-recovery-exact} hinges on the following lemma.  Let $S$ be a superset of the support of $v$.  Relative to the candidate $x = e_1$, any nonzero $\xtilde$ gives components on $\Sc$ that can only increase $\|Vx\|_1$.  Nonzero $\xtilde$ can give components on $S$ that decrease $\|Vx\|_1$.  If the $\ell_1$ norm of $\Vtilde \xtilde$ on $\Sc$ is large enough and the $\ell_1$ norm of $\Vtilde \xtilde$ on $S$ is small enough, then the minimizer  to \eqref{optimization-problem-vx} must be $e_1$.  

\begin{lemma} \label{lemma-exact-main}  Let $V = [v, \Vtilde]$ with  $\|v\|_\infty = 1$, $V_{\istar,:} = [1, \atilde^t]$, supp$(v) \subseteq S$, and $|S| = s$. 
Suppose that  $\| \Vts \xtilde \|_1 \leq 2 s \|\xtilde\|_1$ and $\|\Vtsc \xtilde \|_1 \geq (2 \|\atilde\|_\infty + 2) s \|\xtilde\|_1$ for all $\xtilde$.  Then, $e_1$ is the unique solution to \eqref{optimization-problem-vx}. 
\end{lemma}

\begin{proof}
For any $x$, observe that
\begin{align}
\|Vx\|_1 &= \|v\  x(1) + \Vts \xtilde \|_1 + \| \Vtsc \xtilde \|_1 \\
&\geq    \|v\|_1 |x(1)| - 2 s \| \xtilde \|_1  + \| \Vtsc \xtilde \|_1 \\
&\geq \|v\|_1 |x(1)| + 2 \|\atilde\|_\infty s \|\xtilde \|_1 \label{lower-bound-vx}
\end{align}
where the first inequality is from the upper bound on $\|\Vts \xtilde\|_1$ and the second inequality is from the lower bound on $\|\Vtsc \xtilde\|_1$.  Note that $x = e_1$ is feasible and has value $\| V e_1 \|_1 = \| v\|_1$.  Hence, at a minimizer $\xtildesharp$, 
\begin{align}
\|v\|_1 |\xsharp(1)| + 2 \|\atilde\|_\infty s \|\xtildesharp \|_1 &\leq \|v\|_1.
\end{align}
Using the constraint $\xsharp(1) + \atilde^t \xtildesharp = 1$, a minimizer must satisfy
\begin{align}
\|v\|_1 (1 - \|\atilde\|_\infty \|\xtildesharp\|_1) + 2 \|\atilde\|_\infty s \|\xtildesharp \|_1 &\leq \|v\|_1.
\end{align}
Noting that $\|v\|_1 \leq s$, a minimizer must satisfy 
\begin{align}
2  \|\atilde\|_\infty s \| \xtildesharp \|_1 &\leq \|\atilde\|_\infty s \|\xtildesharp\|_1.
\end{align}
Hence, $\xtildesharp = 0$.  The constraint provides $\xsharp(1) = 1$, which proves that $e_1$ is the unique solution to \eqref{optimization-problem-vx}.  
\end{proof}

To prove Theorem \ref{thm-recovery-exact} by applying Lemma \ref{lemma-exact-main}, we need to study the minimum value of $\|\Vtsc \xtilde\|_1 / \|\xtilde\|_1$ for matrices $\Vtsc$ with i.i.d. ${N}(0,1)$ entries.  Precisely, we will show the following lemma.
\begin{lemma} \label{lemma-min-l1-l1-random}
Let $A$ be a $n \times k$ matrix with i.i.d. ${N}(0,1)$ entries, with $k  \leq n/16$.  There is a universal constant $\ctilde$, such that with high probability, $ \|A x\|_1 /\|x\|_1 \geq \ctilde n/ \sqrt{k}$ for all $x\neq 0$.  This probability is at least $1 - 2 e^{-n/32} - \gamma_1 e^{-\gamma_2 n}$. 
\end{lemma}  

\begin{proof}[Proof of Lemma \ref{lemma-min-l1-l1-random}]
We are to study the problem 
\begin{align}
\min \| A x \|_1 \text{ such that } \|x\|_1 = 1, \label{min-ax-l1}
\end{align}
which is equivalent to 
\begin{align}
\min \| A x \|_1 \text{ such that } \|x\|_1 \geq 1. \label{min-ax-l1-geq}
\end{align}
The minimum value of \eqref{min-ax-l1-geq} can be bounded from below by that of
\begin{align} 
\min \|A x \|_1 \text{ such that } \|x \|_2 \geq 1/\sqrt{k}, \label{min-ax-l2-geq}
\end{align}
because the feasible set of \eqref{min-ax-l1-geq} is included in the feasible set of \eqref{min-ax-l2-geq}.  We now write both the objective and constraint in terms of $Ax$.  To that end, we apply the lower bound in \eqref{sigma-min-random-matrix} to get 
\begin{align}
\P( \|x\|_2 \leq 2  \frac{\|A x \|_2}{ \sqrt{n}}  \text{ for all } x) \geq 1 - 2 e^{-n/32}.
\end{align}
The feasible set of \eqref{min-ax-l2-geq} is contained by the set $\{x \mid \|Ax\|_2 \geq \frac{1}{2} \sqrt{\frac{n}{k}}\}$ with high probability.   Hence, a lower bound to \eqref{min-ax-l2-geq} is with high probability given by
\begin{align}
\min \|Ax \|_1 \text{ such that } \|A x \|_2 \geq \frac{1}{2}\sqrt{\frac{n}{k}}. \label{min-ax-l1-l2}
\end{align}
In order to find a lower bound on \eqref{min-ax-l1-l2}, we apply Lemma \ref{lemma-l1-l2-random-subspace} to the range of $A$, which is a $k$-dimensional random subspace of $\R^n$ with $k \leq n/16$.  Taking $\eta = 1/16$, we see that with high probability, the minimal value of \eqref{min-ax-l1-l2} is bounded from below by $\frac{c_\eta}{2} \frac{n}{\sqrt{k}}$.
The minimal value of \eqref{min-ax-l1-l2}, and hence of \eqref{min-ax-l1}, is bounded from below by $\ctilde n / \sqrt{k}$ for some universal constant $\ctilde$ with probability at least $1 - 2e^{-n/32} - \gamma_1 e^{-\gamma_2 n}$. 
\end{proof}

To prove the theorem by applying Lemma \ref{lemma-exact-main}, we also need to study the maximum value of $\|\Vts \xtilde\|_1 / \|\xtilde\|_1$ for matrices $\Vts$ with i.i.d. ${N}(0,1)$ entries.

\begin{lemma} \label{lemma-max-l1-l1-random}
Let $A$ be a $s \times k$ matrix with i.i.d. ${N}(0,1)$ entries.  Then $\sup_{x\neq 0} \|Ax\|_1 / \|x \|_1 \leq 2s$ with probability at least $1 - k e^{-s}$.
\end{lemma}
\begin{proof}
Note that elementary matrix theory gives that the $\ell_1 \to \ell_1$ operator norm of $A$ is
\begin{align}
\max_{x \neq 0} \frac{\|A x \|_1}{\|x\|_1} = \max_{1 \leq i \leq k} \|A e_i\|_1.
\end{align}
As $Ae_i$ is an $s \times 1$ vector of i.i.d. standard normals, we have
\begin{align}
\P( \|A e_i\|_1 > t ) \leq 2^{s} e^{-t^2/2 s}.
\end{align}
Hence,
\begin{align}
\P(\max_i \|A e_i\|_1 > t ) \leq k 2^{s} e^{-t^2/2s}.
\end{align}
Taking $t = 2 s$, we conclude
\begin{align}
\P(\max_i \|A e_i\|_1 > 2 s ) \leq k 2^{s} e^{-2s} \leq k e^{-s}.
\end{align}

\end{proof}

We can now combine Lemmas \ref{lemma-exact-main}, \ref{lemma-min-l1-l1-random}, and \ref{lemma-max-l1-l1-random} to prove Theorem \ref{thm-recovery-exact}.

\begin{proof}[Proof of Theorem \ref{thm-recovery-exact}]
Let $\ctilde$ be the universal constant given by Lemma \ref{lemma-min-l1-l1-random} and let $c = \ctilde/5$.  We will show that for $\|v\|_0 \leq c \frac{n / \sqrt{\log n}}{\sqrt{k}}$, the minimizer to \eqref{optimization-problem-z-in-v} is $v$ with at least the stated probability.

Let $S$ be any superset of supp$(v)$ with cardinality $s =  \bigl \lfloor c \frac{n / \sqrt{\log n}}{\sqrt{k }} \bigr \rfloor$.  As per Lemma \ref{lemma-exact-main}, $e_1$ is the solution to \eqref{optimization-problem-vx}, and hence $v$ is the unique solution to \eqref{optimization-problem-z-in-v}, if  the following events occur simultaneously:
\begin{align}
\|\Vts \xtilde \|_1 &\leq 2 s \| \xtilde \|_1 \text{ for all } \xtilde \label{proof-condition-vts} \\
\| \atilde \|_\infty &\leq 2 \sqrt{\log n} \label{proof-condition-at-infty}\\
\|\Vtsc \xtilde\|_1 &\geq 5 \sqrt{\log n} s \| \xtilde \|_1 \text{ for all } \xtilde \label{proof-condition-vtsc}
\end{align}
Applying Lemma \ref{lemma-max-l1-l1-random} to the $s \times k$ matrix $\Vts$, we get that \eqref{proof-condition-vts} holds with probability at least \new{$1 - k e^{-s} \geq 1 - k e^{-\lfloor c \sqrt{n/\log n}\rfloor}$.} Classical results on the maximum of a gaussian vector establishes that \eqref{proof-condition-at-infty} holds with probability at least $1 - k / n^2$.  
Because $s \leq n/2$ and $k \leq n/32$, we have that $\Vtsc$ has height at least $n/2$ and width at most $n/32$.  Hence, Lemma \ref{lemma-min-l1-l1-random} gives that $\|\Vtsc \xtilde \|_1 / \|\xtilde\|_1 \geq \ctilde n/\sqrt{k}$ for all $\xtilde\neq 0$ with probability at least $1 - 2 e^{-n/64} - \gamma_1 e^{-\gamma_2 n/2}$.  Because $s \leq \frac{\ctilde}{5} \frac{n/ \sqrt{\log n}}{\sqrt{k}}$, we conclude \eqref{proof-condition-vtsc}, allowing us to apply Lemma \ref{lemma-exact-main}.
Hence, successful recovery occurs with probability at least $1 - 2e^{-n/64} - \gamma_1 e^{-\gamma_2 n/2} - k e^{-\lfloor c \sqrt{n/\log n}\rfloor} - k/n^2 \geq 1 - \new{\tilde{\gamma}_1 e^{-\tilde{\gamma}_2 n}} - k e^{-\lfloor c \sqrt{n/\log n}\rfloor} - k/n^2$, \new{for some $\tilde{\gamma}_1, \tilde{\gamma}_2$.  }

\end{proof}
\subsection{Proof of Theorem \ref{thm-recovery-noisy}}

We will prove the following lemma, of which Theorem \ref{thm-recovery-noisy} is a special case.

\begin{lemma} \label{lemma-recovery-noisy}
Let $k \leq n/32$.  There exists universal constants $c,C$ such that for sufficiently large n,  for all $s \leq c \frac{n / \sqrt{\log n}}{\sqrt{k}}$, and for $i=\istar$, any minimizer $\zsharp$ of  \eqref{optimization-problem-z-in-v} satisfies

\begin{align}
\left \|\zsharp - \frac{v}{v(\istar)} \right \|_2 \leq C \frac{\sqrt{n}}{s} \frac{\|v - \vs \|_1}{\|v \|_\infty} \label{lemma-error-bound}
\end{align}
with probability at least $1 -  \new{\tilde{\gamma}_1 e^{-\tilde{\gamma}_2 n}} - k e^{-s} - k/n^2$. 
\end{lemma}
At first glance, this lemma appears to have poor error bounds for large $n$ and poor probabilistic guarantees for small $s$.  On further inspection, the bounds can be improved by simply considering a larger $s$, possibly even larger than the size of the support of $v$.  Larger values of $s$ simultaneously increase the denominator and decrease the $s$-term approximation error in the numerator of \eqref{lemma-error-bound}.  Taking the largest permissible value $s = \lfloor c \frac{n / \sqrt{\log n}}{\sqrt{k}} \rfloor$, we arrive at Theorem \ref{thm-recovery-noisy}.

Lemma \ref{lemma-recovery-noisy} hinges on the following analog of Lemma \ref{lemma-exact-main}.

\begin{lemma} \label{lemma-noisy-main}  
Fix $1\leq s < n$ and $\alpha > 0$.  Let $V = [v, \Vtilde]$ with  $\|v\|_\infty = 1$, $V_{\istar,:} = [1, \atilde^t]$, $\delta = \|v - \vs\|_1$, supp$(v) \subseteq S$, and $|S| = s$.
If  $\| \Vts \xtilde \|_1 \leq 2 s \|\xtilde\|_1$ and $\|\Vtsc \xtilde \|_1 \geq (2 \|\atilde\|_\infty + 2 + \alpha) s \|\xtilde\|_1$ for all $\xtilde \in \R^{k}$, then any $\xsharp$ minimizing \eqref{optimization-problem-vx} satisfies
\begin{align}
|\xsharp_1 - 1| &\leq \frac{2 \delta}{s}, \qquad \text{ and } \qquad
\|\xtildesharp\|_1 \leq \frac{2 \delta}{s (\|\atilde\|_\infty + \alpha)}.
\end{align}

\end{lemma}

\begin{proof}
For any $x$, observe that
\begin{align}
\|Vx\|_1 &= \| v\cdot x(1)  + \Vts \xtilde \|_1 + \| \Vtsc \xtilde \|_1 \\
&\geq    \|v\|_1 |x(1)| - 2 s \| \xtilde \|_1  + \| \Vtsc \xtilde \|_1 \\
&\geq \|v\|_1 |x(1)| + (2 \|\atilde\|_\infty + \alpha) s \|\xtilde \|_1\\
&\geq  (\|\vs \|_1 - \delta) |x(1)| + (2 \|\atilde\|_\infty + \alpha) s \|\xtilde \|_1
\end{align}
where the first inequality is from the upper bound on $\|\Vts \xtilde\|_1$ and the second inequality is from the lower bound on $\|\Vtsc \xtilde\|_1$.  Note that $x = e_1$ is feasible and has value $\| V e_1 \|_1 = \| v \|_1 \leq \|\vs\|_1 + \delta$.  Hence, at a minimizer $\xtildesharp$, 
\begin{align}
(\|\vs\|_1 - \delta) |\xsharp(1)| + (2 \|\atilde\|_\infty + \alpha) s \|\xtildesharp \|_1 &\leq \|\vs\|_1 + \delta.
\end{align}
Using the constraint $\xsharp(1) + \atilde \xtildesharp = 1$, a minimizer must satisfy
\begin{align}
(\|\vs\|_1 - \delta) (1 - \|\atilde\|_\infty \|\xtildesharp\|_1) + (2 \|\atilde\|_\infty + \alpha) s \|\xtildesharp \|_1 &\leq \|\vs\|_1 + \delta.
\end{align}
Noting that $\|\vs\|_1 \leq s$, a minimizer must satisfy 
\begin{align}
\| \xtildesharp \|_1 \leq \frac{2 \delta}{(\|\atilde\|_\infty + \alpha) s}.
\end{align}
Applying the constraint again, we get
\begin{align}
|\xsharp(1) - 1| \leq \frac{2 \delta}{s}.
\end{align}
\end{proof}

We now complete the proof of Theorem \ref{thm-recovery-noisy} by proving Lemma \ref{lemma-recovery-noisy}.

\begin{proof}[Proof of Lemma \ref{lemma-recovery-noisy}]
Let $\ctilde$ be the universal constant given by Lemma \ref{lemma-min-l1-l1-random} and let $c = \ctilde/6$.  We will show that for any $s \leq c \frac{n / \sqrt{\log n}}{\sqrt{k}}$, the minimizer to \eqref{optimization-problem-z-in-v} is near $v$ with at least the stated probability.

Let $S$ be any superset of supp$(\vs)$ with cardinality $s$.  Applying Lemma \ref{lemma-noisy-main} with $\alpha = \sqrt{\log n}$, we observe that a minimizer $\xsharp$ to \eqref{optimization-problem-vx} satisfies $|\xsharp(1) - 1| \leq 2\delta / s$ and $\| \xtildesharp \|_1 \leq 2\delta / ( s \sqrt{\log n})$ if  the following events occur simultaneously:
\begin{align}
\|\Vts \xtilde \|_1 &\leq 2 s \| \xtilde \|_1 \text{ for all } \xtilde \label{proof-condition-vts-noisy} \\
\| \atilde \|_\infty &\leq 2 \sqrt{\log n} \label{proof-condition-at-infty-noisy}\\
\|\Vtsc \xtilde\|_1 &\geq 6 \sqrt{\log n} s \| \xtilde \|_1 \text{ for all } \xtilde \label{proof-condition-vtsc-noisy}
\end{align}
Applying Lemma \ref{lemma-max-l1-l1-random} to the $s \times k$ matrix $\Vts$, we get that \eqref{proof-condition-vts-noisy} holds with probability at least $1 - k e^{-s}$.  Classical results on the maximum of a gaussian vector establishes that \eqref{proof-condition-at-infty-noisy} holds with probability at least $1 - k / n^2$.  
Because $s \leq n/2$ and $k \leq n/32$, we have that $\Vtsc$ has height at least $n/2$ and width at most $n/32$.  Hence, Lemma \ref{lemma-min-l1-l1-random} gives that $\|\Vtsc \xtilde \|_1 / \|\xtilde\|_1 \geq \ctilde n/\sqrt{k}$ for all $\xtilde\neq 0$ with probability at least $1 - 2 e^{-n/64} - \gamma_1 e^{-\gamma_2 n/2}$.  If $s \leq \frac{\ctilde}{6} \frac{n/ \sqrt{\log n}}{\sqrt{k}}$, we conclude \eqref{proof-condition-vtsc-noisy}, allowing us to apply Lemma \ref{lemma-exact-main}.

It remains to show that $V\xsharp$ is near $v$.  Observe that 
\begin{align}
\|V \xsharp - v\|_2 &= \|V \xsharp - V e_1 \|_2\\
&\leq \|v\|_2 |\xsharp(1) - 1| + \|\Vtilde \xtildesharp \|_2 \\
&\leq \|v\|_2 |\xsharp(1) - 1| + \sigma_\text{max}(\Vtilde) \| \xtildesharp\|_2 \\
&\leq \sqrt{n} |\xsharp(1) - 1| + \frac{3}{2} \sqrt{n} \| \xtildesharp\|_1 \\
&\leq \sqrt{n} \frac{2 \delta}{s}  +  \frac{3}{2} \sqrt{n} \frac{2 \delta}{s \sqrt{\log n}}\\
&\leq C \frac{\sqrt{n}}{s}\delta
\end{align}
The the third inequality uses the fact that $\|v\|_\infty = 1$ and $\sigma_\text{max}(\Vtilde) \leq \frac{3}{2} \sqrt{n}$, which occurs with probability at least $1 - 2e^{-n/32}$ due to the upper bound in \eqref{sigma-min-random-matrix}.
  Hence, approximate recovery occurs with probability at least $1 - 2e^{-n/64} - 2e^{-n/32} - \gamma_1 e^{-\gamma_2 n/2} - k e^{-s} - k/n^2 \new{ \geq 1 - \tilde{\gamma}_1 e^{-\tilde{\gamma}_2 n} - k e^{-s} - k/n^2}$, \new{for some $\tilde{\gamma_1}, \tilde{\gamma}_2$. }

\end{proof}

\section{Simulations} \label{sec:simulations}
\new{
We now present computer simulations that demonstrate when solving the $n$ linear programs \eqref{optimization-problem-z-in-v} can recover the approximately sparse vector $v\in \R^n$ from a planted random subspace \eqref{planted-random-problem}.  In order to obtain a single output from the $n$ programs, we select the output that is the `sparsest' in one of several possible senses.  In this section, we study the effect of these different senses.  We also study the best possible recovery performance of these senses by simulating the behavior of purely random subspaces with no planted sparse vector.
}

\begin{figure}[t!]
\centering 
%
%
\begin{psfrags}%
\psfragscanon%
%
\psfrag{s13}[b][b]{\color[rgb]{0,0,0}\setlength{\tabcolsep}{0pt}\begin{tabular}{c}Sparsity ($s$)\end{tabular}}%
\psfrag{s14}[][]{\color[rgb]{0,0,0}\setlength{\tabcolsep}{0pt}\begin{tabular}{c}Subspace dimensionality ($k$)\end{tabular}}%
\psfrag{s15}[][]{\color[rgb]{0,0,0}\setlength{\tabcolsep}{0pt}\begin{tabular}{c}One program -- $i = \istar$\end{tabular}}%
\psfrag{s16}[][]{\color[rgb]{0,0,0}\setlength{\tabcolsep}{0pt}\begin{tabular}{c}(a)\end{tabular}}%
\psfrag{s17}[b][b]{\color[rgb]{0,0,0}\setlength{\tabcolsep}{0pt}\begin{tabular}{c}Sparsity ($s$)\end{tabular}}%
\psfrag{s18}[][]{\color[rgb]{0,0,0}\setlength{\tabcolsep}{0pt}\begin{tabular}{c}Subspace dimensionality ($k$)\end{tabular}}%
\psfrag{s19}[][]{\color[rgb]{0,0,0}\setlength{\tabcolsep}{0pt}\begin{tabular}{c}All n programs -- min $\ell_1/\ell_\infty$\end{tabular}}%
\psfrag{s20}[][]{\color[rgb]{0,0,0}\setlength{\tabcolsep}{0pt}\begin{tabular}{c}(b)\end{tabular}}%
\psfrag{s21}[b][b]{\color[rgb]{0,0,0}\setlength{\tabcolsep}{0pt}\begin{tabular}{c}Sparsity ($s$)\end{tabular}}%
\psfrag{s22}[][]{\color[rgb]{0,0,0}\setlength{\tabcolsep}{0pt}\begin{tabular}{c}All n programs -- min $\ell_1/\ell_2$\end{tabular}}%
\psfrag{s23}[][]{\color[rgb]{0,0,0}\setlength{\tabcolsep}{0pt}\begin{tabular}{c}Subspace dimensionality ($k$)\end{tabular}}%
\psfrag{s24}[][]{\color[rgb]{0,0,0}\setlength{\tabcolsep}{0pt}\begin{tabular}{c}(c)\end{tabular}}%
\psfrag{s25}[b][b]{\color[rgb]{0,0,0}\setlength{\tabcolsep}{0pt}\begin{tabular}{c}Sparsity ($s$)\end{tabular}}%
\psfrag{s26}[][]{\color[rgb]{0,0,0}\setlength{\tabcolsep}{0pt}\begin{tabular}{c}(d)\end{tabular}}%
\psfrag{s27}[][]{\color[rgb]{0,0,0}\setlength{\tabcolsep}{0pt}\begin{tabular}{c}All n programs -- min thresholded-$\ell_0$\end{tabular}}%
\psfrag{s28}[][]{\color[rgb]{0,0,0}\setlength{\tabcolsep}{0pt}\begin{tabular}{c}Subspace dimensionality ($k$)\end{tabular}}%
%
\psfrag{x01}[t][t]{10}%
\psfrag{x02}[t][t]{20}%
\psfrag{x03}[t][t]{30}%
\psfrag{x04}[t][t]{40}%
\psfrag{x05}[t][t]{50}%
\psfrag{x06}[t][t]{60}%
\psfrag{x07}[t][t]{70}%
\psfrag{x08}[t][t]{80}%
\psfrag{x09}[t][t]{90}%
\psfrag{x10}[t][t]{10}%
\psfrag{x11}[t][t]{20}%
\psfrag{x12}[t][t]{30}%
\psfrag{x13}[t][t]{40}%
\psfrag{x14}[t][t]{50}%
\psfrag{x15}[t][t]{60}%
\psfrag{x16}[t][t]{70}%
\psfrag{x17}[t][t]{80}%
\psfrag{x18}[t][t]{90}%
\psfrag{x19}[t][t]{10}%
\psfrag{x20}[t][t]{20}%
\psfrag{x21}[t][t]{30}%
\psfrag{x22}[t][t]{40}%
\psfrag{x23}[t][t]{50}%
\psfrag{x24}[t][t]{60}%
\psfrag{x25}[t][t]{70}%
\psfrag{x26}[t][t]{80}%
\psfrag{x27}[t][t]{90}%
\psfrag{x28}[t][t]{10}%
\psfrag{x29}[t][t]{20}%
\psfrag{x30}[t][t]{30}%
\psfrag{x31}[t][t]{40}%
\psfrag{x32}[t][t]{50}%
\psfrag{x33}[t][t]{60}%
\psfrag{x34}[t][t]{70}%
\psfrag{x35}[t][t]{80}%
\psfrag{x36}[t][t]{90}%
%
\psfrag{v01}[r][r]{90}%
\psfrag{v02}[r][r]{80}%
\psfrag{v03}[r][r]{70}%
\psfrag{v04}[r][r]{60}%
\psfrag{v05}[r][r]{50}%
\psfrag{v06}[r][r]{40}%
\psfrag{v07}[r][r]{30}%
\psfrag{v08}[r][r]{20}%
\psfrag{v09}[r][r]{10}%
\psfrag{v10}[r][r]{90}%
\psfrag{v11}[r][r]{80}%
\psfrag{v12}[r][r]{70}%
\psfrag{v13}[r][r]{60}%
\psfrag{v14}[r][r]{50}%
\psfrag{v15}[r][r]{40}%
\psfrag{v16}[r][r]{30}%
\psfrag{v17}[r][r]{20}%
\psfrag{v18}[r][r]{10}%
\psfrag{v19}[r][r]{90}%
\psfrag{v20}[r][r]{80}%
\psfrag{v21}[r][r]{70}%
\psfrag{v22}[r][r]{60}%
\psfrag{v23}[r][r]{50}%
\psfrag{v24}[r][r]{40}%
\psfrag{v25}[r][r]{30}%
\psfrag{v26}[r][r]{20}%
\psfrag{v27}[r][r]{10}%
\psfrag{v28}[r][r]{90}%
\psfrag{v29}[r][r]{80}%
\psfrag{v30}[r][r]{70}%
\psfrag{v31}[r][r]{60}%
\psfrag{v32}[r][r]{50}%
\psfrag{v33}[r][r]{40}%
\psfrag{v34}[r][r]{30}%
\psfrag{v35}[r][r]{20}%
\psfrag{v36}[r][r]{10}%
%
\resizebox{14cm}{!}{\includegraphics{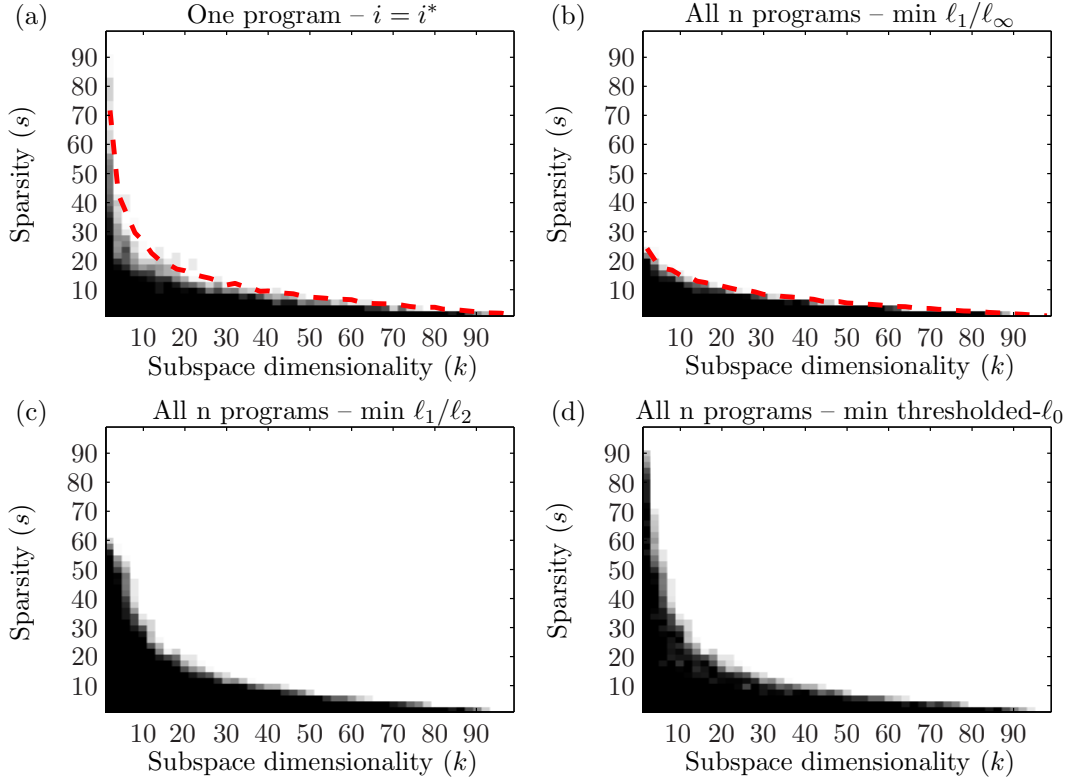}}%
\end{psfrags}%
%

  \caption{ \setdoublespace 
Empirical probability of recovery versus approximate sparsity and subspace dimensionality $k$.  For a given $s$, the vector to be recovered is a noisy version of $1_S$, where $|S|=s$.  For all values of $k$ and $s$,  we solve $n$ programs of form \eqref{optimization-problem-z-in-v}, one for each $1\leq i \leq n$.  Panel (a) shows the output corresponding to $i = \istar$.  Panels (b), (c), and (d) show the output of the $n$ programs that was smallest in the $\ell_1/\ell_\infty$, $\ell_1/\ell_2$, and thresholded-$\ell_0$ senses, respectively.  The diagram shows the probability of recovery, as measured by \new{50} independent trials.  White represents recovery with probability zero.  Black represents recovery with probability 1.  The dashed lines in panels (a) and (b) represent upper bounds for the maximal recoverable sparsity based on the behavior of random subspaces with no planted vector.}
    \label{fig:phase-transition}
\end{figure}

The parameters for the linear programs \eqref{optimization-problem-z-in-v} are as follows.  Let $n=100$,  $1 \leq s \leq n$, and $S = \{1, \ldots, s\}$.  Let $1_S$ be the vector that is $1$ on $S$ and $0$ on $S^c$.    
 We attempt to recover the approximately $s$-sparse vector $v = 1_S + \new{\delta u} $, where \new{$\delta=0.01$} and $u$ has i.i.d. Gaussian entries and is normalized such that $\|u\|_1 = 1$.   Let $\istar = \argmax_i |v(i)|$.   We solve \eqref{optimization-problem-z-in-v} for $1 \leq i \leq n$ using YALMIP \cite{L2004} with the SDPT3 solver \cite{TTT1999, TTT2003}.   Among these $n$ outputs, we let $\zsharp$ be the one that (a) corresponds to $i = \istar$; (b) is smallest in the sense of $\ell_1/\ell_\infty$; (c) is smallest in the sense of $\ell_1/\ell_2$; or (d) is sparsest in the sense of thresholded-$\ell_0$ at the level \new{$\epsilon = 0.01$}.  We will refer to (a) as the oracle selector because it corresponds to an oracle that tells us the index of the largest component of $v$.    We call a recovery successful if $\left \|\zsharp - v / v(\istar) \right \|_2 \leq \new{0.01}$.  Figure \ref{fig:phase-transition} shows the probability of successful recovery, as computed 
over \new{50} independent trials, for  many values of  $k$ and the approximate sparsity $s$.    
Near and below the phase transitions, simulations were performed for all even values of $k$ and $s$.  In the large regions to the top-right of the phase transitions, simulations were performed only for values of $k$ and $s$ that are multiples of 5.  In this region, the probably of recovery was always zero.  

\new{We also ran computer simulations in the noiseless case $\delta=0$.  The outputs of these simulations are not shown because they look the same as Figure \ref{fig:phase-transition} for the parameters as above. 
}

\new{To get an indication for when \eqref{optimization-problem-z-in-v} can successfully recover an approximately sparse vector from a planted random subspace, we also study the behavior of the corresponding problems without a planted vector.}  Specifically, we compute
\begin{align}
\min \|z\|_1 \text{ such that } z \in \text{span} \{\vtilde_{1}, \ldots, \vtilde_{k}\}, z(i) = 1, \label{min-in-random-space}
\end{align}
for i.i.d. $\vtilde_{j} \sim {N}(0, I_n)$ and for all $1\leq i \leq n$.  The curve in Figure \ref{fig:phase-transition}a shows the dependence on $k$ of the minimal value of \eqref{min-in-random-space} for a single fixed $i$, as found by the median over 50 independent trials.  The curve in Figure \ref{fig:phase-transition}b shows the solution to \eqref{min-in-random-space} that is smallest in the  $\ell_1/\ell_\infty$ sense, also found by the median of  50 trials.  We note that the least $\ell_1/\ell_\infty$ solution to \eqref{min-in-random-space} is also a solution of the problem
\begin{align}
\min \frac{\|z\|_1}{\|z\|_\infty} \text{ such that } z \in \text{span} \{\vtilde_{1}, \ldots, \vtilde_{k}\}, z\neq 0. \label{min-in-random-space-l1linf}
\end{align}
\new{These curves represent upper bounds for the sparsity of a recoverable signal: for sparsities of $v$ above this curve, there will be linear combinations of random vectors that are considered sparser than $v$ in the provided sense.  Hence, recovering $v$ would be impossible.    These curves are also upper bounds in the noiseless case of $\delta = 0$. }  Figures \ref{fig:phase-transition}c and \ref{fig:phase-transition}d do not have corresponding lines because we do not have a tractable method for directly optimizing sparsity in the $\ell_1/\ell_2$ or thresholded-$\ell_0$ senses.

Our primary observation is that the selection process for the sparsest output among the $n$ programs \eqref{optimization-problem-z-in-v} can greatly effect recovery performance, especially for planted-random subspaces with low dimensionality.  Figures \ref{fig:phase-transition}a and \ref{fig:phase-transition}b show that $v$ can be recovered by the $i = \istar$ program yet discarded as not sparse enough in the $\ell_1/\ell_\infty$ sense.   \new{In the parameter regime we simulated, the $\ell_1/\ell_\infty$ selector exhibits the worst overall recovery performance.  The thresholded-$\ell_0$ selector gives the best performance, though it has the drawback of having a threshold parameter.}

Our second observation is that recovery can succeed even when the $i=\istar$ instance of \eqref{optimization-problem-z-in-v} fails. That is, solving all $n$ programs of form \eqref{optimization-problem-z-in-v} can outperform the result of a single program, even when an oracle tells us the index of the largest coefficient.  \new{ For example, compare the phase transitions for small $k$ of Figure \ref{fig:phase-transition}a to Figures \ref{fig:phase-transition}c and \ref{fig:phase-transition}d.  This effect is more likely when $k$ is small and $v$ has many large components.  To see why, note that successful recovery is expected when $\|\Vtilde_{i,:}\|_\infty$ is small.  If $k$ is small, $\|\Vtilde_{i,:}\|_\infty$ is more likely to be small.  If there are many indices $i$ where $v(i)$ is large, it is likely that $\| \Vtilde_{i,:}\|_\infty$ will be small enough for successful recovery with some other value of $i$.   }

\new{
Our third observation is that purely random subspaces provide good estimates for the best recoverable sparsity with a planted random subspace, under the oracle and $\ell_1/\ell_\infty$ selectors.  This agreement is revealed by the similarity of the phase transition and dashed lines in Figures \ref{fig:phase-transition}a and \ref{fig:phase-transition}b.  In the plots, the region of successful recovery on average can not exceed the dashed lines because recovery requires that the planted sparse vector be smaller in the appropriate sense than any vector in the purely random part of the subspace.  Because the dashed lines are also upper bounds in the noiseless case of $\delta = 0$, the performance of the recovery process can not be significantly better than that shown in Figures \ref{fig:phase-transition}a and \ref{fig:phase-transition}b.  
}

Potentially, the recovery region could be improved significantly beyond that indicated by Figure \ref{fig:phase-transition}.  For example, one could consider  more than $n$ linear programs, each with a normalization against a different random direction in $\R^n$.  Such an approach immediately gives rise to a natural tradeoff: recovery performance  may be improved at the expense of  more linear programs that need to be solved.  We leave this relationship for future study.

\section*{Funding}
This work was supported by the National Science Foundation to P.H. and L.D.;  the Alfred P. Sloan Foundation to L.D.; TOTAL S.A. to L.D.; the Air Force Office of Scientific Research to L.D.; and the Office of Naval Research to L.D.  
\section*{Acknowledgements} The authors would like to thank Jonathan Kelner and Vladislav Voroninski for helpful discussions.

\bibliographystyle{plain}
\bibliography{subspace-refs}

\clearpage


\end{document}